\documentclass[english,conference, a4paper]{IEEEtran}
\usepackage[T1]{fontenc}
\usepackage[latin9]{inputenc}
\usepackage{prettyref}
\usepackage{float}
\usepackage{amsthm}
\usepackage{amsmath}
\usepackage{amssymb}
\usepackage{graphicx}
\usepackage{esint}
\usepackage{xargs}[2008/03/08]

\makeatletter

\floatstyle{ruled}
\newfloat{algorithm}{tbp}{loa}
\providecommand{\algorithmname}{Algorithm}
\floatname{algorithm}{\protect\algorithmname}

\theoremstyle{plain}
\newtheorem{thm}{\protect\theoremname}
\theoremstyle{plain}
\newtheorem{prop}[thm]{\protect\propositionname}
\theoremstyle{remark}
\newtheorem{rem}[thm]{\protect\remarkname}

\IEEEoverridecommandlockouts

\usepackage{dimadefs}

\newcommand{\RR}{\ensuremath{\mathbb{R}}}

\newcommand{\nparams}{\ensuremath{R}}

\let\oldhat\hat
\renewcommand{\vec}[1]{\mathbf{#1}}
\renewcommand{\hat}[1]{\oldhat{\mathbf{#1}}}

\usepackage{prettyref}
\newrefformat{sub}{Subsection \ref{#1}}
\newrefformat{subfig}{Subfigure \ref{#1}}
\newrefformat{app}{Appendix \ref{#1}}
\newrefformat{cor}{Corollary \ref{#1}}
\newrefformat{alg}{Algorithm \ref{#1}}
\newrefformat{step}{step \ref{#1}}
\newrefformat{def}{Definition \ref{#1}}
\newrefformat{clm}{Claim \ref{#1}}
\newrefformat{prop}{Proposition \ref{#1}}
\newrefformat{prob}{Problem \ref{#1}}
\newrefformat{que}{Question \ref{#1}}
\newrefformat{rem}{Remark \ref{#1}}

\@ifundefined{showcaptionsetup}{}{%
 \PassOptionsToPackage{caption=false}{subfig}}
\usepackage{subfig}
\makeatother

\usepackage{babel}
\providecommand{\propositionname}{Proposition}
\providecommand{\remarkname}{Remark}
\providecommand{\theoremname}{Theorem}

\begin{document}

\title{Algebraic signal sampling, Gibbs phenomenon and Prony-type systems}

\author{\IEEEauthorblockN{Dmitry Batenkov\IEEEauthorrefmark{1}\IEEEauthorrefmark{2}
and Yosef Yomdin\IEEEauthorrefmark{1}\IEEEauthorrefmark{3}}\IEEEauthorblockA{\IEEEauthorrefmark{1}Department
of Mathematics,Weizmann Institute of Science, Rehovot 76100, Israel}\IEEEauthorblockA{\IEEEauthorrefmark{2}Email:
dima.batenkov@weizmann.ac.il}\IEEEauthorblockA{\IEEEauthorrefmark{3}Email:
yosef.yomdin@weizmann.ac.il}\thanks{This research was supported
by the Adams Fellowship Program of the Israeli Academy of Sciences
and Humanities, ISF Grant No. 639/09 and by the Minerva foundation.}}
\maketitle
\begin{abstract}
Systems of Prony type appear in various signal reconstruction problems
such as finite rate of innovation, superresolution and Fourier inversion
of piecewise smooth functions. We propose a novel approach for solving
Prony-type systems, which requires sampling the signal at arithmetic
progressions. By keeping the number of equations small and fixed,
we demonstrate that such ``decimation'' can lead to practical improvements
in the reconstruction accuracy. As an application, we provide a solution
to the so-called Eckhoff's conjecture, which asked for reconstructing
jump positions and magnitudes of a piecewise-smooth function from
its Fourier coefficients with maximal possible asymptotic accuracy
-- thus eliminating the Gibbs phenomenon.
\end{abstract}
\global\long\def\np{\ensuremath{K}}
\global\long\def\jp{\ensuremath{x}}
\global\long\def\jc{\ensuremath{a}}
\newcommandx\fc[1][usedefault, addprefix=\global, 1=k]{\ensuremath{c_{#1}}}
\global\long\def\sc{\ensuremath{M}}
\global\long\def\fun{\ensuremath{f}}
\newcommandx\er[1][usedefault, addprefix=\global, 1=k]{\delta_{#1}}
\global\long\def\err{\varepsilon}
\newcommandx\meas[1][usedefault, addprefix=\global, 1=k]{m_{#1}}
\newcommandx\nnmeas[1][usedefault, addprefix=\global, 1=k]{\widetilde{m}_{#1}}
\global\long\def\nmeas{S}
\newcommandx\apprmeas[1][usedefault, addprefix=\global, 1=k]{\widetilde{m}_{#1}}
\newcommandx\frsum[2][usedefault, addprefix=\global, 1=\fun, 2=\sc]{\mathfrak{F}_{#2}\left(#1\right)}
\global\long\def\smooth{\ensuremath{\Psi}}
\global\long\def\sing{\ensuremath{\Phi}}
\global\long\def\nn#1{\widetilde{#1}}
\global\long\def\ord{\ensuremath{d}}
\global\long\def\jcc{\alpha}
\global\long\def\scc{N}
\global\long\def\w{\ensuremath{\omega}}

\global\long\def\noims{\tilde{\vec{y}}}
\global\long\def\exams{\vec{y}}
\global\long\def\src{\vec{x}}
\global\long\def\noisrc{\tilde{\src}}
\global\long\def\ncoeffs{\mathcal{C}}
\global\long\def\df{p}
\global\long\def\init{t}
\newcommandx\fwm[1][usedefault, addprefix=\global, 1=I]{\mathcal{P}_{#1}}

\section{Introduction}

The ``Prony system'' of equations
\begin{equation}
\meas=\sum_{j=1}^{\np}c_{j}z_{j}^{k},\qquad c_{j},z_{j}\in\complexfield,\; k\in\naturals\label{eq:basic-prony}
\end{equation}
appeared originally in the work of R.Prony \cite{prony1795essai}
in the context of fitting a sum of exponentials to observed data samples.
He showed that the unknowns $\left\{ c_{j,}z_{j}\right\} _{j=1}^{\np}$
can be recovered explicitly from $\left\{ \meas[0],\dots,\meas[2\np-1]\right\} $
by what is known today as ``Prony's method''. The system \eqref{eq:basic-prony}
appears in areas such as frequency estimation, Padé approximation,
array processing, statistics, interpolation, quadrature, radar signal
detection, error correction codes, and many more. In modern signal
processing, \eqref{eq:basic-prony} is of fundamental importance in
the field of sub-Nyquist sampling (related terms are superresolution
\cite{candes2012towards,donoho1992superresolution} and finite rate
of innovation \cite{dragotti2007sma}). A basic problem there is to
recover an unknown ``spike train'', a linear combination of $\delta$-functions
\[
\fun\left(x\right)=\sum_{j=1}^{\np}b_{j}\delta\left(x-x_{j}\right),\qquad c_{j}\in\reals,\; x_{j}\in\left[-\pi,\pi\right]
\]
from its Fourier samples
\begin{equation}
\widehat{f}\left(k\right)=\frac{1}{2\pi}\int_{-\pi}^{\pi}\fun\left(t\right)\ee^{-\imath kt}\dd t.\label{eq:fourier-coeffs}
\end{equation}

The resulting system is of course a special case of \eqref{eq:basic-prony}.
If a more general model is considered,
\begin{align}
f(x) & =\sum_{j=1}^{\np}\sum_{\ell=0}^{\ell_{j}-1}b_{\ell,j}\der{\delta}{\ell}(x-\jp_{j}),\quad b_{\ell,j}\in\reals,\;\jp_{j}\in\left[-\pi,\pi\right],\label{eq:gen-delta-fun}
\end{align}
then \eqref{eq:fourier-coeffs} becomes, after a change of variables,
\begin{equation}
\meas=\sum_{j=1}^{\np}z_{j}^{k}\sum_{\ell=0}^{\ell_{j}-1}c_{\ell,j}k^{\ell},\qquad c_{\ell,j}\in\complexfield,\;\left|z_{j}\right|=1.\label{eq:polynomial-prony}
\end{equation}
Many research efforts are devoted to stable solution of Prony-type
systems (see e.g. \cite{badeau2008performance,batenkov2011accuracy,donoho2006stable,potts2010parameter,stoica1989music}
and references therein). We propose a novel approach to this problem,
which requires sampling the signal at arithmetic progressions. By
keeping the number of equations small and fixed, we demonstrate (in
\prettyref{sec:dec-prony}) that such ``decimation'' can lead to
practical improvements in the reconstruction accuracy, to a certain
extent avoiding a well-known numerical instability of these systems.

In \prettyref{sec:fourier-1d} we consider the problem of recovering
a piecewise-smooth function, including the positions of its discontinuities,
from its Fourier samples. The algebraic reconstruction method due
to K.Eckhoff in essense required a solution of a particular instance
of the system \eqref{eq:polynomial-prony} with the error in the left-hand
side having a certain asymptotic decay rate. Previously it was shown
in \cite{Batenkov2011,batyomAlgFourier} that this approach yields
a nonlinear approximation which is ``half as accurate'' compared
to the best possible bound. As we elaborate in \prettyref{sec:fourier-1d},
applying the decimation technique to the Prony-type system results
in full asymptotic accuracy, thus completely eliminating the Gibbs
phenomenon.

In \prettyref{sec:future-work} we discuss several promising directions
for future research.

\section{\label{sec:dec-prony}Decimated Prony-type systems}

Suppose that the ``polynomial Prony model'' \eqref{eq:polynomial-prony}
is to be fitted to the noisy measurements $\nnmeas[0],\dots\nnmeas[\sc-1]$.
We denote the number of unknowns by $\nparams=\sum_{j=1}^{\np}\left(\ell_{j}+1\right)$.
At first sight, using all the $\sc$ measurements for fitting should
improve reconstruction accuracy. While this is certainly justified
in the case where the noise statistics are known (as demonstrated
in e.g. \cite{badeau2008performance,stoica1989music}), this might
backfire if the noise is ``adversary'', or ``worst-case''. Potts
\& Tasche \cite{potts2010parameter} show that when Prony system \eqref{eq:basic-prony}
is solved by least squares minimization for all $\sc$ equations at
once, then even if the nodes $\left\{ z_{j}\right\} $ are detected
very accurately, the error for magnitudes is amplified by a factor
of $\sqrt{\sc\nparams}$. This shows that it might actually be productive
to stay with small number of measurements. We are therefore justified
in making a simplifying assumption that the number of equations used
for reconstruction equals the number of unknowns $\nparams$. In this
case the solution to the reconstruction problem can be characterized
as the exact inversion of the measurement mapping $\fwm:\complexfield^{\nparams}\to\complexfield^{\nparams}$
which associates to any parameter vector $\src=\bigl\{\{c_{ij}\},\{\jp_{i}\}\bigr\}\in\complexfield^{\nparams}$
its corresponding exact measurement vector $\exams=\left(\meas[i_{0}],\dots,\meas[i_{\nparams-1}]\right)\in\complexfield^{\nparams}$
where $I=\left\{ i_{0}<i_{1}<\dots<i_{\nparams-1}\right\} \subset\left[0,\sc-1\right]$
is a given index set. Perhaps the most natural choice for the index
sets $I$ is given by arithmetic progressions
\[
I_{\init,\df}=\left\{ \init,\init+\df,\dots,\init+\left(\nparams-1\right)\df\right\} ,\qquad\init\geqslant0,\;\df\geqslant1.
\]

Following \cite{batenkov2011accuracy}, we estimate for such $I=I_{\init,\df}$
the (local) stability of inversion by the \emph{Lipschitz constant
of $\fwm^{-1}$} at the regular points of $\fwm$, which in turn are
given by the following proposition.
\begin{prop}
\label{prop:regular-pts}The vector $\vec{x}=\left(\left\{ z_{j},c_{i,j}\right\} \right)\in\complexfield^{\nparams}$
is a regular point of $\fwm$ with $I=I_{\init,\df}$ if and only
if $z_{j}^{\df}\neq z_{i}^{\df}$ for $i\neq j$, and $c_{\ell_{j}-1,j}\neq0$
for all $j=1,\dots,\np$.
\end{prop}
We have the following upper bound on the accuracy of \emph{any solution
method}.
\begin{thm}
\label{thm:decimated-prony-accuracy}Consider the polynomial Prony
system \eqref{eq:polynomial-prony} with a fixed structure $\left\{ \np,\left\{ \ell_{j}\right\} _{j=1}^{\np}\right\} $
on $I=I_{\init,\df}$, and let $\vec{x}=\left(\left\{ z_{j},c_{i,j}\right\} \right)\in\complexfield^{\nparams}$
be a regular point of $\fwm[I]$. If the error in each measurement
is bounded in absolute value by $\err\ll1$, then the errors in recovering
the components of the original parameter vector $\vec{x}$ satisify{\footnotesize{
\begin{eqnarray*}
\left|\Delta c_{i,j}\right| & \leq & C\left(i,\ell_{j}\right)\left(\frac{2}{\delta_{\df}}\right)^{\nparams}\left(\frac{1}{2}+\frac{R}{\delta_{\df}}\right)^{\ell_{j}}\frac{\init^{\ell_{j}-i}}{\df^{i}}\left(1+\frac{\left|c_{i-1,j}\right|}{\left|c_{\ell_{j}-1,j}\right|}\right)\err,\\
\left|\Delta z_{j}\right| & \leq & \frac{2}{\ell_{j}!}\left(\frac{2}{\delta_{\df}}\right)^{\nparams}\frac{1}{\left|c_{\ell_{j}-1,j}\right|}\df^{-\ell_{j}}\err,
\end{eqnarray*}
}}where $\delta_{\df}\isdef\min_{i\neq j}\left|z_{j}^{\df}-z_{i}^{\df}\right|$
and $C\left(i,\ell_{j}\right)$ is an explicit constant (for consistency
we take $c_{-1,j}=0$ in the above formula).
\end{thm}
This result directly generalizes earlier stability estimates of \cite{batenkov2011accuracy}
for the special case $I=I_{0,1}$. The proofs of both \prettyref{prop:regular-pts}
and \prettyref{thm:decimated-prony-accuracy} are based on factorizing
the Jacobian matrix of the map $\fwm$ along the same lines as in
\cite{batenkov2011accuracy}, while adding the analysis of the Jacobian's
dependence on $\init$ and $\df$.

\begin{figure}
\subfloat[ESPRIT, $\ord=3$]{\includegraphics[width=0.9\columnwidth]{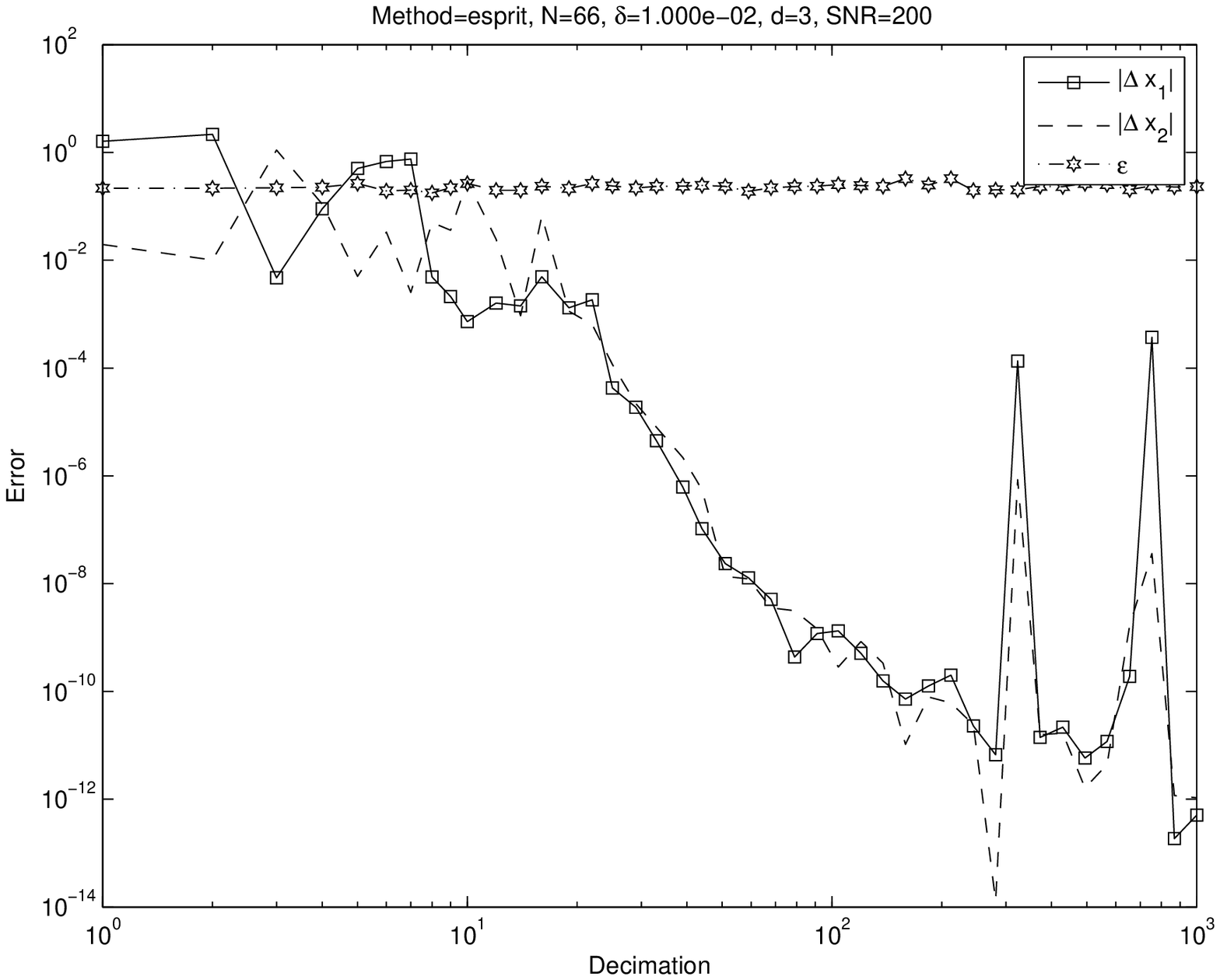}

}

\subfloat[LS, $\ord=3$]{\includegraphics[width=0.9\columnwidth]{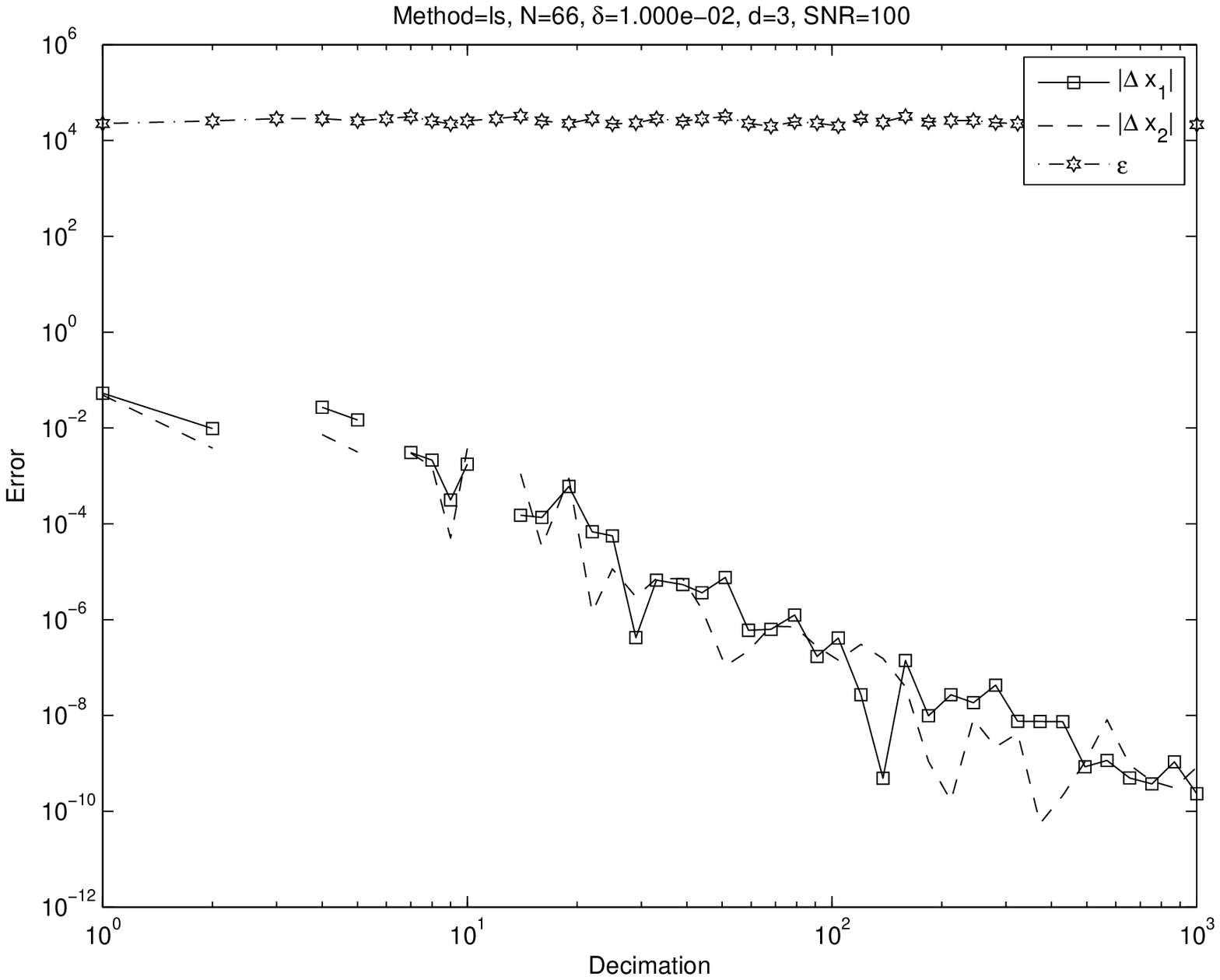}

}

\caption{Reconstruction error as a function of the decimation with fixed number
of measurements ($\sc=66$). The signal has two nodes with distance
$\delta=10^{-2}$ between each other. Notice that ESPRIT requires
significantly higher Signal-to-Noise Ratio in order to achieve the
same performance as LS.}
\label{fig:fixed-number-meas}

\end{figure}

\begin{figure}
\subfloat[ESPRIT, $\ord=3$]{\includegraphics[width=0.9\columnwidth]{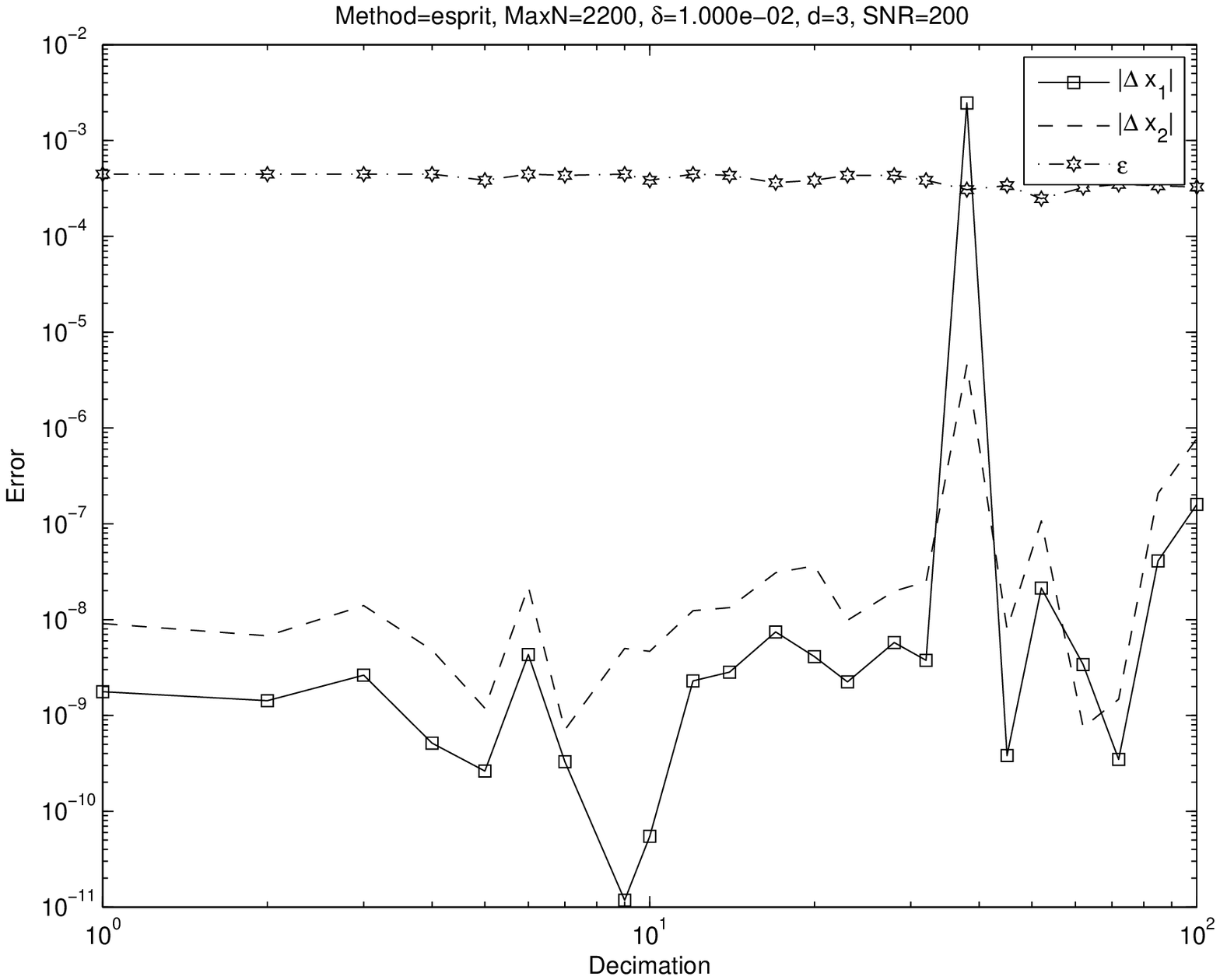}

}

\subfloat[LS, $\ord=3$]{\includegraphics[width=0.9\columnwidth]{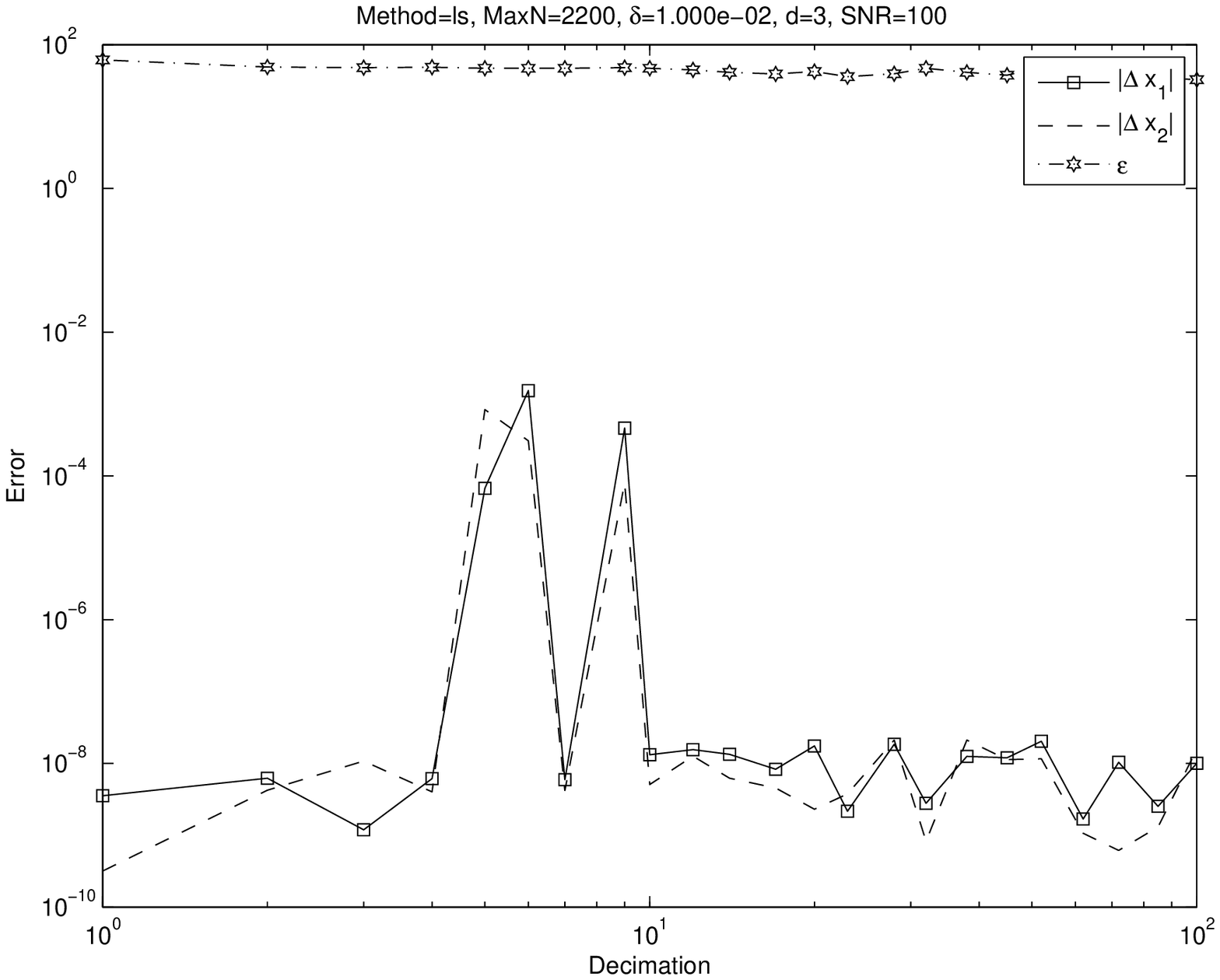}}

\caption{Reconstruction error as a function of the decimation, reducing number
of measurements from $\sc=2200$ to $\sc=22$. The signal has two
nodes with distance $\delta=10^{-2}$ between each other. The reconstruction
accuracy remains almost constant.}
\label{fig:full-decimation}
\end{figure}

Now suppose that the number of available measurements $\sc\to\infty$,
while the noise $\err$ remains bounded. It is easy to see that for
the index set $I=I_{0,\left\lfloor \frac{\sc}{\nparams}\right\rfloor }$
we obtain an improvement in accuracy of recovering the jump $z_{j}$
of the order $\sim\sc^{\ell_{j}}$, compared with the non-decimated
measurement set $I_{0,1}$.
\begin{rem}
If initially two nodes are close (say by $\delta$), the decimation
improves accuracy up to a certain limit. To see this, just substitute
$\delta_{\df}\sim\df\delta$ into \prettyref{thm:decimated-prony-accuracy}
and get an improvement by factor of $\df^{-\nparams-\ell_{j}}$. 
\end{rem}
Turning to particular solution methods, the decimation is fairly straightforward
to implement. Indeed, taking any algorithm for the standard Prony-type
system, one just needs to make the following modifications (for simplicity
we consider only the recovery of the nodes $\left\{ z_{j}\right\} $).
\begin{enumerate}
\item Choose the decimation parameter $\df$.
\item Feed the original algorithm with the decimated measurements $\meas[0],\meas[\df],\meas[2\df]\dots,$
and obtain the estimated node $w_{j}$.
\item Take $z_{j}=\sqrt[\df]{w_{j}}$.
\end{enumerate}
We have tested the decimation technique according to the above procedure
on two well-known algorithms for Prony systems - ESPRIT \cite{badeau2008performance}
and nonlinear least squares (LS, implemented by MATLAB's \texttt{lsqnonlin}).
In the first experiment, we fixed the number of measurements to be
66, and changed the decimation parameter $\df$, while keeping the
noise level constant. The accuracy of recovery increased with $\df$
-- see \prettyref{fig:fixed-number-meas}. In the second experiment,
we fixed the highest available measurement to be $\sc=2200$, and
changed the decimation from $\df=1$ to $\df=100$ (thereby reducing
the number of measurements from $2200$ to just $22$). The accuracy
of recovery stayed relatively constant -- see \prettyref{fig:full-decimation}.
Note that such a reduction leads to a corresponding decrease in the
running time (calculating singular-value decomposition of large matrices,
as in ESPRIT, is a time-consuming operation).

\section{\label{sec:fourier-1d}Piecewise-smooth Fourier reconstruction}

Consider the problem of reconstructing an integrable function $\fun:\left[-\pi,\pi\right]\to\RR$
from a finite number of its Fourier coefficients \eqref{eq:fourier-coeffs}.
If $f$ is $C^{d}$ and periodic, then the truncated Fourier series
$\frsum\isdef\sum_{|k|=0}^{\sc}\fc(\fun)\ee^{\imath kx}$ approximates
$f$ with error at most $C\cdot\sc^{-d-1}$, which is optimal. If,
however, $f$ is not smooth even at a single point, the rate of accuracy
drops to only $\sc^{-1}$. Still, one can hope to restore the best
accuracy by using the a-priori information to produce some non-standard
summation method. This accuracy problem, also known as the Gibbs phenomenon,
is very important in applications, such as calculation of shock waves
in  PDEs. It has received much attention especially in the last few
decades - see e.g. a recent book \cite{jerriGibbs11}.

The so-called ``algebraic approach'' to this problem, first suggested
by K.Eckhoff \cite{eckhoff1995arf}, is as follows. Assume that $\fun$
has $\np>0$ jump discontinuities $\left\{ \jp_{j}\right\} _{j=1}^{\np}$
, and $\fun\in C^{\ord}$ in every segment $\left(\jp_{j-1},\jp_{j}\right)$.
We say that in this case $\fun$ belongs to the class $PC\left(\ord,\np\right)$.
Denote the associated jump magnitudes at $\jp_{j}$ by $\jc_{\ell,j}\isdef\der{\fun}{\ell}(\jp_{j}^{+})-\der{\fun}{\ell}(\jp_{j}^{-}).$
Then write the piecewise smooth $\fun$ as the sum $\fun=\smooth+\sing$,
where $\smooth(x)$ is smooth and periodic and $\sing(x)$ is a piecewise
polynomial of degree $\ord$, uniquely determined by $\left\{ \jp_{j}\right\} ,\left\{ \jc_{\ell,j}\right\} $
such that it ``absorbs'' all the discontinuities of $\fun$ and its
first $\ord$ derivatives. In particular, the Fourier coefficients
of $\sing$ have the explicit form{\small{
\begin{equation}
\fc(\sing)=\frac{1}{2\pi}\sum_{j=1}^{\np}\ee^{-\imath k\jp_{j}}\sum_{\ell=0}^{\ord}(\imath k)^{-\ell-1}\jc_{\ell,j},\quad k=1,2,\dots.\label{eq:singular-fourier-explicit}
\end{equation}
}}For $k\gg1$, we have $\left|\fc\left(\sing\right)\right|\sim k^{-1}$,
while $\left|\fc\left(\smooth\right)\right|\sim k^{-\ord-2}$. Consequently,
Eckhoff suggested to pick large enough $k$ and solve the approximate
system of equations \eqref{eq:polynomial-prony} where $\meas=2\pi\left(\imath k\right)^{\ord+1}\fc\left(\fun\right)$,
$z_{j}=\ee^{-\imath\jp_{j}}$ and $c_{\ell,j}=\imath^{\ell}\jc_{\ord-\ell,j}$.
His proposed method of solution was to use the known values $\left\{ \meas\right\} _{k\in I}$
where
\begin{equation}
I=\left\{ \sc-\left(\ord+1\right)\np+1,\sc-\left(\ord+1\right)\np+2,\dots,\sc\right\} \label{eq:eckhoff-index-set}
\end{equation}
to construct an algebraic equation satisfied by the unknowns $\left\{ z_{1},\dots,z_{\np}\right\} $,
and solve this equation numerically. Based on some explicit computations
for $\ord=1,2;\;\np=1$ and large number of numerical experiments,
he conjectured that his method would reconstruct the jump locations
with accuracy $\sc^{-\ord-1}$.

Let us consider the problem in the framework of Prony type system
\eqref{eq:polynomial-prony}. The error term is of magnitude $\left|\err\right|\sim\sc^{-1}$.
The index set \eqref{eq:eckhoff-index-set} is just $I_{\init,\df}$
with $\init\sim\sc,\;\df=1$ (i.e. no decimation). Therefore, by \prettyref{thm:decimated-prony-accuracy}
we get accuracy only of order $\left|\Delta\jp_{j}\right|\sim\sc^{-1}$. 

Now consider the decimated setting for this problem. By the above,
we can approximate eash jump $\jp_{j}$ up to accuracy $\sc^{-1}$.
Set
\[
\scc=\left\lfloor \frac{\sc}{\left(\ord+2\right)\np}\right\rfloor .
\]
Now take the index set $I_{\init,\df}$ where $\init=\df=\scc$, i.e.
$I_{\scc,\scc}=\left\{ \scc,2\scc,\dots,\sc\right\} .$ As before,
$\left|\epsilon\right|\sim\sc^{-1}$, but now due to decimation we
get accuracy $\left|\Delta\jp_{j}\right|\sim N^{-\ord-1}N^{-1}\sim\sc^{-\ord-2}.$
In \cite{batFullFourier,batyomAlgFourier} we develop an algorithm
(see \prettyref{alg:full-fourier}) which in fact attains this accuracy.
This result can be summarized as follows.

\begin{algorithm}
Let $f\in PC\left(\ord,\np\right)$, and assume that $\fun=\sing^{\left(\ord\right)}+\smooth$
where $\sing^{\left(\ord\right)}$ is the piecewise polynomial absorbing
all discontinuities of $\fun$, and $\smooth\in C^{\ord}.$
\begin{enumerate}
\item Obtain initial approximations for $\left\{ \jp_{1},\dots,\jp_{\np}\right\} $
by any standard method (i.e. Eckhoff's method of order zero).
\item Localize each $\jp_{j}$ by multiplying with a mollifier (convolution
in Fourier domain).
\item Solve resulting Prony system with $\np=1$ and $\init=\df=\left\lfloor \frac{\sc}{d+2}\right\rfloor $
(decimation).
\item Take the final approximation to be{\footnotesize{
\[
\begin{split}\widetilde{\fun}= & \nn{\Phi}\left(\left\{ \nn{\jc}_{\ell,j},\nn{\jp}_{j}\right\} \right)\\
 & +\sum_{\left|k\right|\leq\sc}\left\{ \fc(\fun)-\frac{1}{2\pi}\sum_{j=1}^{\np}\ee^{-\imath\nn{\jp_{j}}k}\sum_{\ell=0}^{d}\frac{\widetilde{\jc}_{\ell,j}}{(\imath k)^{\ell+1}}\right\} \ee^{\imath kx}.
\end{split}
\]
}}{\footnotesize \par}
\end{enumerate}
\caption{Full accuracy Fourier reconstruction of piecewise smooth functions}
\label{alg:full-fourier}
\end{algorithm}

\begin{thm}
Let $f\in PC\left(\ord,\np\right)$, so that $\fun=\sing^{\left(\ord\right)}+\smooth$
where $\sing^{\left(\ord\right)}$ is the piecewise polynomial with
Fourier coefficients \eqref{eq:singular-fourier-explicit}, and $\smooth\in C^{\ord}.$
Assume that there exist constants $J,A,B,R$ such that
\begin{eqnarray*}
\min_{i\neq j}\left|\jp_{i}-\jp_{j}\right| & \geq & J>0,\quad\left|\fc\left(\smooth\right)\right|\leq R\cdot k^{-\ord-2},\\
\left|\jc_{\ell,j}\right| & \leq & A<\infty,\quad\left|\jc_{0,j}\right|\geq B>0.
\end{eqnarray*}
Then the approximation $\widetilde{\fun}$ obtained by \prettyref{alg:full-fourier}
satisfies for $\sc\gg1${\small{
\[
\begin{split}\left|\nn{\jp_{j}}-\jp_{j}\right| & \leq C_{1}\left(\ord,\np,J,A,B,R\right)\cdot\sc^{-\ord-2};\\
\left|\nn{\jc}_{\ell,j}-\jc_{\ell,j}\right| & \leq C_{2}\left(\ord,\np,J,A,B,R\right)\cdot\sc^{\ell-\ord-1},\quad0\leqslant\ell\leqslant\ord;\\
\left|\nn{\fun}\left(x\right)-\fun\left(x\right)\right| & \leq C_{3}\left(\ord,\np,J,A,B,R\right)\cdot\sc^{-\ord-1}.
\end{split}
\]
}}{\small \par}
\end{thm}
Note that the pointwise bound $\left|\fun\left(x\right)-\nn{\fun}\left(x\right)\right|$
is valid ``away from discontinuities''. Some numerical experiments,
elaborated in \cite{batFullFourier,batyomAlgFourier}, confirm these
theoretical accuracy predictions.

\section{\label{sec:future-work}Future work}

Stable solution of Prony-type systems in the most general setting
must take into account the possibility of colliding nodes. We believe
that a reparametrization of the equations in the basis of finite differences
is a promising approach to this problem. We have obtained initial
results in \cite{byPronySing12,yom2009Singularities}, and plan to
continue in this direction.

The Fourier inversion problem for piecewise-analytic functions is
still widely open (see e.g. \cite{adcock2012stability}). While our
results provide spectral convergence in this setting, it is still
unknown if the algebraic method can be pushed to exponential or at
least root-exponential accuracy.

Edge detection from spectral data is a well-researched problem, see
e.g. \cite{engelberg2008recovery,greengard2011spectral} and references
therein. We expect that the 1D procedure can be generalized to treat
the general case via some form of a ``separation'', or ``slice
reconstruction'' (see e.g. \cite{BatGolYom2011} for an example of
such a method, dealing with reconstruction from moments).

\bibliographystyle{plain}
\bibliography{../../../bibliography/all-bib}

\end{document}